\documentclass[11pt]{amsart}
\usepackage{pstricks}
\usepackage{latexsym}
\usepackage{amssymb}

 \textwidth=\textwidth \textheight=\textheight
\theoremstyle{plain}
\newtheorem{thm}{Theorem}
\newtheorem{theorem}[thm]{Theorem}
\newtheorem{lem}[thm]{Lemma}
\newtheorem{prop}[thm]{Proposition}

\newtheorem{lemma}[thm]{Lemma}

\theoremstyle{definition}
\newtheorem{definition}[thm]{Definition}

\newcommand{\set}[1]{ \left\{ #1 \right\} }

\def\linv{{\mathrm{linv}}}
\def\R{{\mathcal R}}

\def\ll{\lambda}

\def\d{\delta}

\def\zz{{\mathbb Z}}

\def\uqsln{U_q(\widehat{ \mathfrak{
sl}}_n)}

\def\g{{\mathcal G}}

\def\r{{\mathrm r}}

\def\sh{{\rm sh}}
\def\Rib{{\rm Rib}}
\def\inv{{\rm inv}}
\def\quot{{\rm quot}}

\def\diag{{\rm diag}}

\newcommand{\be}[1]{\begin{equation} \label{#1}}
\newcommand{\ee}{\end{equation}}

\newcommand{\spin}{\mathrm{spin}}

\newcommand{\wt}{\text{wt}}

\psset{unit=1pt} \psset{arrowsize=4pt 1} \psset{linewidth=1pt}
\newrgbcolor{mygreen}{.2 .7 .2}
\newrgbcolor{myred}{.7 .3 .2}
\newgray{mygray}{.90}
\countdef\x=23 \countdef\y=24 \countdef\z=25 \countdef\t=26

\def\abox(#1,#2)#3{
\x=#1 \y=#2 \multiply\x by 16 \multiply\y by 16 \z=\x \t=\y
\advance\z by 16 \advance\t by 16
\psline(\x,\y)(\x,\t)(\z,\t)(\z,\y)(\x,\y) \advance\x by 8
\advance\y by 8 \rput(\x,\y){{\bf #3}}}

\def\hdom(#1,#2)#3{
\x=#1 \y=#2 \multiply\x by 16 \multiply\y by 16 \z=\x \t=\y
\advance\z by 32 \advance\t by 16
\psline(\x,\y)(\x,\t)(\z,\t)(\z,\y)(\x,\y) \advance\x by 16
\advance\y by 8 \rput(\x,\y){{\bf #3}}}

\def\vdom(#1,#2)#3{
\x=#1 \y=#2 \multiply\x by 16 \multiply\y by 16 \z=\x \t=\y
\advance\z by 16 \advance\t by 32
\psline(\x,\y)(\x,\t)(\z,\t)(\z,\y)(\x,\y) \advance\x by 8
\advance\y by 16 \rput(\x,\y){{\bf #3}}}

\def\rec(#1,#2,#3,#4){
\psline(#1,#2)(#3,#2)(#3,#4)(#1,#4)(#1,#2) }

\begin{document}
\date{July~11, 2004}
\title{Two results on domino and ribbon tableaux}
\author{Thomas Lam}
\address{Department of Mathematics,
         M.I.T., Cambridge, MA 02139}

\email{thomasl@math.mit.edu}
\begin{abstract}
Inspired by the spin-inversion statistic of Schilling, Shimozono
and White \cite{SSW} and Haglund et al. \cite{HHLRU} we relate the
symmetry of ribbon functions to a result of van Leeuwen, and also
describe the multiplication of a domino function by a Schur
function.
\end{abstract}
\maketitle
\section{Introduction}
Lascoux, Leclerc and Thibon \cite{LLT} defined spin-weight
generating functions $\g_{\ll/\mu}^{(n)}(X;q)$ (from hereon called
\emph{ribbon functions}) for ribbon tableaux.  They proved that
these functions were symmetric functions using the action of the
Heisenberg algebra on the Fock space of $\uqsln$.  For the $n=2$
case of domino tableaux, a combinatorial proof of the symmetry and
in fact a description of the expansion of
$\g_{\ll/\mu}^{(n)}(X;q)$ in terms of Schur functions is given by
the Yamanouchi domino tableaux of Carr\'{e} and Leclerc \cite{CL}.
More recently, Schilling, Shimozono and White \cite{SSW} and
separately Haglund et. al. \cite{HHLRU} have described the spin
statistic of a ribbon tableau in terms of an inversion number on
the $n$-quotient. This article gives two applications of this
inversion number towards the ribbon functions.

Our first application is a proof of the symmetry of ribbon
functions using a result of van Leeuwen \cite{vL1} developed from
his spin-preserving Knuth correspondence for ribbon tableaux.  The
result says roughly that the spin generating functions for adding
horizontal ribbon strips above or below a lattice path vertical on
both ends are equal.  Another `elementary' but more systematic
proof of the symmetry of ribbon functions will appear in
\cite{Lam}.

Our second application is an imitation of Stembridge's concise
proof of the Littlewood Richardson rule \cite{Ste} for the domino
tableau case.  We describe the expansion of $s_\nu(X)
\g_{\mu/\rho}^{(2)}(X;q)$ in the basis of Schur functions in terms
of \emph{$\nu$-Yamanouchi domino tableaux}. This description
appears to be new and also gives a  shorter proof of the result of
Carr\'{e} and Leclerc \cite{CL}, corresponding to $\nu = (0)$, the
empty partition .

In the last section we describe explicitly a bijection in terms of
words required to prove the symmetry of ribbon functions.

\medskip

We refer the reader to \cite{LLT,Lam1} for the necessary
definitions and notation concerning ribbon tableaux, spin and
ribbon functions.  We will always think of our partitions and
tableaux as being drawn in the English notation.

{\bf Acknowledgements.} This project is part of my Ph.D. Thesis
written under the guidance of Richard Stanley.  I am grateful for
all his advice and support over the last couple of years.

\section{Spin-inversion statistic}
We will use the spin-inversion statistic from \cite{HHLRU} as its
description is considerably shorter than the one in \cite{SSW},
and we will only be interested in how spin changes rather than its
exact value. Let $\quot_n(T) = (T^{(0)},\ldots,T^{(n-1)})$ denote
the $n$-quotient of a ribbon tableau $T$ (which may have skew
shape). With the $n$-core fixed, semistandard ribbon tableaux are
in bijection with such $n$-tuples of usual tableaux.  The diagonal
$\diag(s)$ of a cell $s \in \quot_n(T)$ is equal to the diagonal
of $T$ on which the head of the corresponding ribbon $\Rib(s)$
lies.  For a cell $s \in T^{(i)}$ it is given by $\diag(s) = nc(s)
+ c_i$ for some offsets $c_i$ depending on the $n$-core of
$\sh(T)$.  Here $c(s) = j - i$ is the usual \emph{content} of a
square $s = (i,j)$.  An \emph{inversion} is a pair of entries
$T(x) = a, T(y) = b$ such that $a<b$ and $0 < \diag(x) - \diag(y)
< n$. We denote by $\inv(T) = \inv(\quot_n(T))$ the number of
inversions of $\quot_n(T)$. We have \cite{HHLRU}
\begin{lemma}
\label{lem:HHLRU} Given a skew shape $\ll/\mu$, there is a
constant $e(\ll/\mu)$ such that for every standard $n$-ribbon
tableau $T$ of shape $\ll/\mu$, we have $\spin(T) = e(\ll/\mu) -
\inv(\quot_n(T))$.
\end{lemma}

We shall use a particular \emph{diagonal reading order} on our
tableaux.  Let $T$ be a ribbon tableaux.  The \emph{reading word}
$\r(T)$ is given by reading the diagonals of $\quot_n(T)$ in
descending order, where in each diagonal the larger numbers are
read first. We will regularly abuse notation by allowing ourselves
to identify ribbons in $T$ with squares of the $n$-quotient
$\quot_n(T)$.  We will also identify a skew shape $\ll/\mu$ which
is a horizontal ribbon strip with the correpsonding horizontal
ribbon strip tableau $T$ satisfying $\sh(T) = \ll/\mu$.

\section{Symmetry of ribbon functions}
We fix the length $n \geq 1$ of our ribbons throughout.

Recall that the standard way to prove that a Schur function is
symmetric is to give an involution $\alpha_i$ on semistandard
tableaux of shape $\ll$ which swaps the number of $i$'s and
$(i+1)$'s, for each $i$.  This is known as the Bender-Knuth
involution. Our first aim is to study the symmetry of the ribbon
functions $\g_{\ll/\mu}^{(n)}(X;q)$ from the perspective of the
$n$-quotient. This symmetry is equivalent to the existence of a
\emph{ribbon Bender-Knuth involution} $\sigma_i$ on ribbon
tableaux $T$ which changes the number of $i$'s and $i+1$'s while
preserving spin.

We call a skew shape $\ll/\mu$ a double horizontal ribbon strip if
it can be tiled by two horizontal ribbon strips.  Let
$\R^{a,b}_{\ll/\mu}$ be the set of ribbon tableaux of shape
$\ll/\mu$ filled with $a$ $1$'s and $b$ $2$'s.  To obtain a ribbon
Bender-Knuth involution, it suffices to find a spin preserving
bijection between $\R^{a,b}_{\ll/\mu}$ and $\R^{b,a}_{\ll/\mu}$
for every $a$ and $b$ and every double horizontal strip $\ll/\mu$.
Let $T \in \R^{a,b}_{\ll/\mu}$. Suppose some tableau $T^{(i)}$ of
the $n$-quotient contains a column with two squares, then those
two squares must be $1$ on top of a $2$.

We first show that we may reduce to the case that the $n$-quotient
contains no such columns.  If $(x,y)$ is an inversion of $T$ we
say that the inversion \emph{involves} $x$ and $y$.  Let
$\inv_x(T)$ denote the number of inversions of $T$ which involve
$x$.

\begin{lem}
\label{lem:onetwo} Let $T$ be a ribbon tableau and $\quot_n(T)$
contain two squares $x$ and $y$ in the same column such that $T(x)
= i$ and $T(y) = i+1$.  Let $T'$ be a semistandard ribbon tableau
obtained from $T$ by changing a `$i$' to a `$i+1$'. Then
\[\inv_x(T) + \inv_y(T) = \inv_x(T') + \inv_y(T').\]
\end{lem}
\begin{proof}
We first note that $\diag(x) = \diag(y) + n$.  Thus the only
relevant inversions come from squares $z$ satisfying $\diag(x) >
\diag(z) > \diag(y)$ and $T(z) \in \set{i,i+1}$.  We check
directly that regardless of the value of $T(z)$, the cell $z$
contributes exactly one inversion to $\inv_x(T) + \inv_y(T)$ and
thus to $\inv_x(T') + \inv_y(T')$ as well.
\end{proof}

Lemma \ref{lem:onetwo} combined with Lemma \ref{lem:HHLRU} shows
that to prove that all ribbon functions are symmetric functions we
only need to check it for horizontal ribbon strips $\ll/\mu$.  For
a horizontal ribbon strip $\ll/\mu$, let $I_{\ll/\mu} \subset \zz$
be the set of diagonals such that $\quot_n(\ll/\mu)$ contains a
cell. It follows from Lemma \ref{lem:HHLRU} that the symmetry of
$\g_{\ll/\mu}^{(n)}(X;q)$ implies the symmetry for all horizontal
strips $\nu/\rho$ with the same set of diagonals $I_{\nu/\rho} =
I_{\ll/\mu}$ -- only the constant $e(\nu/\rho)$ has changed. It is
easy to see that given a set of diagonals $I \subset \zz$, we can
find a horizontal ribbon strip $\ll/\mu$ such that $I_{\ll/\mu} =
I$ and so that $\ll/\mu$ is tileable using vertical ribbons only.
Thus the symmetry of all ribbon functions reduces to the symmetry
of ribbon functions $\g^{(n)}_{\ll/\mu}(X;q)$ corresponding to a
horizontal ribbon strip $\ll/\mu$ tileable only using vertical
ribbons.  In fact it is clear that we need only check this
symmetry for such shapes which are connected.

\section{Connection with a result of van Leeuwen}
Curiously, the symmetry of these special ribbon functions follows
from a result of van Leeuwen concerning adding ribbons above and
below a fixed lattice path. We identify the steps of an infinite
lattice path $P$ going up and right with a doubly infinite
sequence $p = \set{p_i}_{i=-\infty}^{\infty}$ of $0$'s and $1$'s,
where a $0$ corresponds to a step to the right and a $1$
corresponds to a step up.  We may think of such lattice paths as
the boundary of a shape (or partition) in which case the bit
string is known as the edge sequence \cite{vL2}.  For our
purposes, the indexing of $\set{p_i}$ is unimportant.

Van Leeuwen's result is the following \cite[Claim 1.1.1]{vL1}.
\begin{prop}
\label{prop:vL} Let $p = \set{p_i}_{i=-\infty}^{\infty}$ be a
lattice path which is vertical at both ends. Let $R_p$ denote the
generating function
\[
R_p(X,q) = \sum_S q^{\spin(S)}X^{|S|}
\]
where the sum is over all horizontal ribbon strips $S$ that can be
attached below $p$.  Let $\tilde{p}$ denote $p$ reversed.  Then
\[
R_p(X,q) = R_{\tilde{p}}(X,q).
\]
\end{prop}
Note that the generating functions $R_p(X,q)$ are finite, since
only finitely many horizontal ribbon strips can be placed under a
lattice path which is vertical at both ends.  The lattice path
$\tilde{p}$ should be thought of as rotating $p$ upside-down, so
that $R_{\tilde{p}}(X,q)$ enumerates the ways of adding a
horizontal ribbon strip above $p$ (see \cite{vL1}).

We will also need the following technical lemma \cite[Lemma
5.2.2]{vL1} to make a calculation with spin.  For a set $I \subset
\zz$  of diagonals, we denote $\spin_I(T)$ to be the sum of the
spins of the ribbons of $T$ whose heads lie on the diagonals of
$I$.
\begin{lemma}[\cite{vL1}]
\label{lem:vL}  Let $\ll$, $\mu$, $\nu$ be partitions so that
$\ll/\mu$, $\ll/\nu$, $\mu/\nu$ are all horizontal ribbon strips.
Let $I, J \subset \zz$ be the set of diagonals occurring in
$\ll/\mu$ and $\mu/\nu$ respectively.  Then
\[
\spin_I(\ll/\nu) - \spin(\ll/\mu) = \spin_J(\ll/\nu) -
\spin(\mu/\nu).
\]
\end{lemma}

\begin{prop}
\label{prop:sym} Let $\ll/\nu$ be a connected skew shape which is
tileable by vertical ribbons only.  Then
$\g_{\ll/\nu}^{(n)}(x_1,x_2;q)$ is a symmetric function.
\end{prop}
\begin{proof}
In the notation of Proposition \ref{prop:vL}, we pick $p$ so that
$\ll/\nu$ is the shape obtained by adding as many vertical ribbons
as possible below $p$ to give a horizontal ribbon strip.
Alternatively, we can think of $\ll/\nu$ as the bounded region
obtained by shifting the lattice path upwards $n$ steps.  Let $m =
|\ll/\nu|/n$. Let $S_1$ be a horizontal ribbon strip with $a \leq
m$ ribbons added below $p$ which we assume has shape $\mu/\nu$.
Filling $S_1$ with 1's there is a unique way to add another
horizontal ribbon strip $S_2$ filled with 2's to give a tableau $T
\in \R^{a,b}_{\ll/\nu}$.

Since $\spin_I(\ll/\nu) = (n-1)|I|$ for any valid set of diagonals
$I \subset I_{\ll/\nu}$, we have $\spin(S_2) = (n-1)(2a - m) +
\spin(S_1)$ by Lemma \ref{lem:vL}.  Summing over all $S_1$, we get
\[\g_{\ll/\nu}^{(n)}(x_1,x_2;q) = x_2^m q^{-(n-1)m}
R_p\left(\frac{x_1}{x_2}q^{2(n-1)},q^2\right).\] However, we can
also obtain the tableau $T$ by counting the horizontal ribbon
strip $S_2$ containing $2$ first, so a similar argument gives
$\g_{\ll/\nu}^{(n)}(x_1,x_2;q) = x_1^m q^{-(n-1)m}
R_{\tilde{p}}\left(\frac{x_2}{x_1}q^{2(n-1)},q^2\right)$. Since
$R_p = R_{\tilde{p}}$ by Proposition \ref{prop:vL} we obtain
$\g_{\ll/\nu}^{(n)}(x_1,x_2;q) = \g_{\ll/\nu}^{(n)}(x_2,x_1;q)$.
\end{proof}

The following theorem follows immediately from Proposition
\ref{prop:sym} and earlier discussion.

\begin{theorem}
\label{thm:sym}  Let $\ll/\mu$ be any skew shape tileable by
$n$-ribbons. Then $\g_{\ll/\mu}^{(n)}(X;q)$ is a symmetric
function.
\end{theorem}

Theorem \ref{thm:sym} was first shown by Lascoux, Leclerc and
Thibon \cite{LLT} using an action of the Heisenberg algebra on the
Fock space of $\uqsln$.

\section{Generalised Yamanouchi domino tableaux}
In this section we imitate a proof of the Littlewood Richardson
rule due to Stembridge \cite{Ste}, which we apply to domino
tableaux.  We fix $n=2$ throughout this section. Define the
generalised (domino) $q$-Littlewood Richardson coefficients
$c^{\ll}_{\mu/\rho, \nu}(q)$ by
\[
s_\nu(X) \g_{\mu/\rho}(X;q) = \sum_\ll c^{\ll}_{\mu/\rho, \nu}(q)
s_\ll(X).
\]

Let $\set{\sigma_r}$ denote a set of fixed domino Bender-Knuth
involutions which exist by Theorem \ref{thm:sym}.  Let $w = w_1
w_2 \cdots w_k$ be a sequence of integers. Then the weight $\wt(w)
= (\wt_1(w), \wt_2(w), \ldots)$ is the composition of $k$ such
that $\wt_i(w) = |\set{j \mid w_j = i}|$. If $T$ is a ribbon
tableau, let $T_{\geq j}$ and $T_{> j}$ denote the set of ribbons
lying in diagonals which are $\geq j$ and $> j$ respectively (and
similarly for $T_{< j}$ and $T_{\leq j}$). These are not tableaux,
but the compositions $\wt(T_{\geq j})$ and $\wt(T_{>j})$ are well
defined, in the usual manner.

\begin{definition}
Let $\ll$ be a partition.  A word $w = w_1 w_2 \cdots w_k$ is
$\ll$-Yamanouchi if for any initial string $w_1 w_2 \cdots w_i$,
and any integer $l$, we have $\wt_l(w_1\cdots w_i) + \ll_l \geq
\wt_{l+1}(w_1\cdots w_i) + \ll_{l+1}$.  A domino tableau $D$ is
\emph{$\ll$-Yamanouchi} if its reading word $\r(D)$ is
$\ll$-Yamanouchi.
\end{definition}
One can check that $(0)$-Yamanouchi is essentially the notion of
Yamanouchi introduced by Carr\'{e} and Leclerc \cite{CL}.

\begin{theorem}
The generalised $q$-Littlewood Richardson coefficients are given
by
\[
c^{\ll}_{\mu/\rho, \nu}(q) = \sum_Y q^{\spin(Y)}
\]
where the sum is over all $\nu$-Yamanouchi domino tableaux $Y$ of
shape $\mu/\rho$ and weight $\ll$.
\end{theorem}

\begin{proof}
Our proof will follow Stembridge's proof \cite{Ste} nearly step by
step.  We will prove the Theorem in the variables $x_1,\ldots,x_m$
and will always think of a tableau $D$ in terms of its
$2$-quotient.  By definition,
\[
\g_{\mu/\rho}(X;q) = \sum_D q^{\spin(D)}x^D
\]
where the sum is over all semistandard domino tableaux of shape
$\mu/\rho$ filled with numbers in $[1,m]$.  Let $a_{\ll+\d}$
denote the alternating sum $\sum_w (-1)^w x^{w(\ll+\d)}$ where the
sum is over all permutations $w \in S_m$.  Then
\begin{align}
 a_{\ll+\d}\g_{\mu/\rho}(X;q) &= \sum_w \sum_D
q^{\spin(D)} (-1)^w x^{D+w(\ll+\d)}\\ &= \label{eq:fun} \sum_D
q^{\spin(D)} \sum_w
 (-1)^w x^{w(D+\ll+\d)}\\ &= \sum_D q^{\spin(D)} a_{D+
\ll + \d}. \label{eq:bad}
\end{align}
To obtain (\ref{eq:fun}) we have used Theorem \ref{thm:sym} to see
that the weight generating function for domino tableaux with fixed
spin is $w$ invariant.  We call $D$ a Bad Guy if
\[
\ll_k + \wt_k(D_{> j}) < \ll_{k+1} +\wt_{k+1}(D_{\geq j})
\]
for some $j$ and $k$.  Of all such pairs $(j,k)$, we pick one that
maximises $j$ and amongst those we pick the smallest $k$.  Thus
the reading word of $\r(D_{>j})$ is $\ll$-Yamanouchi and the
$j$-th diagonal of $D$ contains a $k+1$ (and possibly a $k$) while
the $(j+1)$-th diagonal contains no $k$.

Now let $S$ be the set of dominoes obtained from $D_{< j}$ by
including the $k$ on the $j$-th diagonal if any.  Set $S^* =
\sigma_k(S)$. This makes sense since the squares of $S$ containing
a $k$ or $k+1$ forms a double horizontal strip which is actually
of skew shape, so we can apply the Bender-Knuth involution.  Now
since $\sh(S) = \sh(S^*)$ we can attach $S^*$ back onto $D_{\geq
j}$ to obtain a tableau $D^*$. We check that $D^*$ is a
semistandard domino tableau. This is the case as only $k$'s and
$k+1$'s are changed into each other, and the boundary diagonals
$j$ and $j+1$ only contain $k+1$'s (there are two conditions to
check, one for each tableau of the 2-quotient). Also note that if
there is a $k$ in diagonal $j$ of $S$ then there must be a $k+1$
immediately below it, so it will always remain a $k$ in $S^*$.

It follows immediately from the construction that $D \mapsto D^*$
is an involution on the set of Bad Guys.  We check that it is
spin-preserving by counting the number of inversions.  Since we
have assumed that $\sigma_k$ preserves spin, the only inversions
that we have to be concerned about are those where $D(x) = k+1$
and $D(y) = k$ and $\diag(x) = j-1$ and $\diag(y)= j$. But if the
$j$-th diagonal contains a $k$, then there is a $k+1$ immediately
below it, so by Lemma \ref{lem:onetwo}, it can be ignored for
calculations of spin in $D$, $D^*$ and also $S$ and $S^*$.  So
$\spin(D) = \spin(D^*)$.

Now,
\[
a_{D+ \ll+\d} = - a_{D^* + \ll + \d},
\]
since $s_k(D+ \ll+\d) = D^*+\ll+\d$, so the contributions of the
Bad Guys to the sum (\ref{eq:bad}) cancel out.  The tableaux which
are not Bad Guys are exactly the $\ll$-Yamanouchi tableaux.
Dividing both sides of (\ref{eq:bad}) by $a_{\d}$ and using the
bialternant formula $s_\ll(X) = a_{\ll+\d}/a_{\d}$ now gives the
Theorem.
\end{proof}

Unfortunately, this proof seems to fail for ribbon tableaux with
$n > 2$.  The similarly defined involution $T \mapsto T^*$ no
longer preserves either semistandard-ness or spin.

We should remark also that Carr\'{e} and Leclerc's algorithm
mapping a domino tableau $D$ to a pair $(Y,T)$ of a Yamanouchi
domino tableau and a usual Young tableau can also be interpreted
in terms of the $2$-quotient.

\section{Word sequence formulation of ribbon function symmetry}
We end the paper by describing explicitly the bijection needed to
prove symmetry of ribbon functions in terms of certain sequences.
Let $n \geq 1$ be an integer.

\begin{definition}
A \emph{$(1,2,\emptyset)$-word} is a sequence $(a_1,a_2, \ldots,
a_m)$ where each $a_i \in \set{1,2,\emptyset}$, such that whenever
$a_i = 2$, then $a_{i+n} \neq 1$.  The \emph{form} $F_a$ of a
sequence $(a_1,a_2,\ldots,a_m)$ is the finite set $F_a = \set{i
\in [1,m] \mid a_i = \emptyset}$. The weight $\wt(a)$ of such a
word $a = (a_1,\ldots,a_m)$ is $(\mu_1,\mu_2)$ where $\mu_i =
\#\set{j:a_j = i}$.
\end{definition}

\begin{definition}
A $n$-\emph{local inversion} of a $(1,2,\emptyset)$-word
$(a_1,a_2, \ldots, a_m)$ is a pair $(i,j)$ satisfying $1 \leq i <
j \leq m$ and $j-i < n$ such that $a_i = 2$ and $a_j = 1$. We let
$\linv_n(w)$ denote the number of $n$-local inversions of $w$.
\end{definition}

The following proposition makes the connection between
$(1,2,\emptyset)$-words and a ribbon Bender Knuth involution.

\begin{prop}
\label{prop:word} The symmetry of ribbon functions is equivalent
to the following identity on $(1,2,\emptyset)$-words for each
positive integer $m$, form $F \subset [1,m]$ and weight
$(\mu_1,\mu_2)$:
\begin{equation}
\label{eq:main} \sum_{a: \wt(a) = (\mu_1,\mu_2)}q^{\linv_n(a)} =
\sum_{a: \wt(a) = (\mu_2,\mu_1)}q^{\linv_n(a)}
\end{equation}
where the sum is over all $(1,2,\emptyset)$-words with length $m$,
form $F$ and specified weight.
\end{prop}
\begin{proof}
We have already established that we need only be concerned with
tableaux which are horizontal ribbon strips filled with ribbons
labelled $1$ and $2$. Our $(1,2,\emptyset)$-words are simply the
(reversed) reading words of these ribbon tableaux where the form
$F$ keeps track of the empty diagonals. The Proposition follows
immediately from Lemma \ref{lem:HHLRU}.
\end{proof}
We remark that when the form $F$ is the emptyset, a bijection
giving (\ref{eq:main}) is obtained by reversing the sequence and
changing $2$'s to $1$'s and vice versa.

\end{document}